\documentclass[english]{article}
\usepackage[T1]{fontenc}
\usepackage[latin1]{inputenc}
\usepackage{amssymb}

\makeatletter

\providecommand{\tabularnewline}{\\}

\usepackage{babel}
\makeatother
\begin{document}
\noindent \textbf{\Large The Sum and Chain Rules for Maximal Monotone
Operators}{\Large \par}

\strut

\strut

\noindent {\large M.D. VOISEI}{\large \par}

\noindent \emph{The University of Texas -- Pan American, Department
of Mathematics, Edinburg, Texas - 78541, USA}

\strut

\strut

\strut

\noindent \textbf{Abstract.} This paper is primarily concerned with
the problem of maximality for the sum $A+B$ and composition $L^{*}ML$
in non-reflexive Banach space settings under qualifications constraints
involving the domains of $A,B,M$. Here $X$, $Y$ are Banach spaces
with duals $X^{*}$, $Y^{*}$, $A,B:X\rightrightarrows X^{*}$, $M:Y\rightrightarrows Y^{*}$
are multi-valued maximal monotone operators, and $L:X\rightarrow Y$
is linear bounded. Based on the Fitzpatrick function, new characterizations
for the maximality of an operator as well as simpler proofs, improvements
of previously known results, and several new results on the topic
are presented.

\strut


\strut


\eject

\noindent \textbf{1. Introduction}

\strut

\noindent The problem concerning the maximality of the sum of two
maximal monotone operators was first stated and solved by Rockafellar
in reflexive Banach spaces followed by a sum rule for the convex subdifferential
in general Banach spaces (see {[}8, Theorem 1], {[}15, Theorem 2.8.3]).
At the same time the conjecture which states that the reflexivity
of the space can be avoided was formulated. Later, the sum rule for
full-space domain operators was proved by Heisler (see {[}10, Theorem
37.4] or Theorem 4.1 below) and for single-valued linear operators
by Phelps and Simons (see {[}6, Theorem 7.2]). Recently, Voisei {[}11,12,13]
proved similar calculus rules for closed convex domain monotone operators
in non-reflexive Banach spaces under weaker forms of the qualification
constraint and solved completely the linear case (see {[}12,13] or
Corollary 5.14. below). Using a topological argument, the sum rule
for operators with the intersection of their domain interiors non-empty
was shown to hold by Borwein (see {[}1]). Chain rules in the context
of reflexive Banach spaces were obtained by Penot {[}5], Z$\breve{{\rm a}}$linescu
{[}14], and Borwein {[}1].

In the present note we want to shed a new light on the ideas of the
proofs and present new points of view as well as simpler arguments,
improvements of some of the past results, and new results concerning
the maximality of the sum or of the precomposition with a linear operator
in the non-reflexive Banach space setting. 

The plan of the paper is as follows. Next section introduces the Fitzpatrick
and Penot functions together with their main features. In Section
3 new characterizations for the maximality and representability of
an operator are discussed. Section 4 contains a simple proof of Heisler's
result. Section 5 deals with the calculus of maximal monotone operators.
This paper concludes with some improvements of the results contained
in Voisei {[}12,13] and several other new results on the topic.

\strut

\strut

\noindent \textbf{2. The natural dual system. The Fitzpatrick and
Penot functions}

\strut

Let $X$ be a Banach space with dual $X^{*}$ and bi-dual $X^{**}$.
A multi-valued operator $A:D(A)\subset X\rightrightarrows X^{*}$
is called \emph{monotone} if, for every $x_{1},x_{2}\in D(A)$\emph{,}
$x_{1}^{*}\in Ax_{1}$, $x_{2}^{*}\in Ax_{2}$\begin{equation}
\langle x_{1}-x_{2},x_{1}^{*}-x_{2}^{*}\rangle\ge0,\label{eq:}\end{equation}
where\[
p(x,x^{*})=\langle x,x^{*}\rangle=x^{*}(x),\ (x,x^{*})\in X\times X^{*},\]
stands for the duality pairing in $X\times X^{*}$.

For the sake of notation simplicity we identify operators with their
graphs and write\begin{equation}
x\in D(A),\ x^{*}\in Ax\Leftrightarrow(x,x^{*})\in A.\label{eq:}\end{equation}
With this notation $A$ is monotone iff $p(a_{1}-a_{2})\ge0$ for
every $a_{1},a_{2}\in A$.

A monotone operator is considered \emph{maximal monotone} if it is
maximal in the sense of inclusion in $X\times X^{*}$.

Let $Z=X\times X^{*}$. The natural dual system is formed by $(Z,Z)$
with the dual product\begin{equation}
(x,x^{*})\cdot(y,y^{*})=x^{*}(y)+y^{*}(x),\ x,y\in X,\ x^{*},y^{*}\in X^{*}.\label{eq:}\end{equation}

The convex conjugate with respect to the natural duality of $f:Z\rightarrow\mathbb{R}\cup\{\infty\}$
is given by\begin{equation}
f^{*}(w)=\sup\{ w\cdot z-f(z);\ z\in Z\},\ w\in Z,\label{eq:}\end{equation}
Notice that \begin{equation}
z\cdot z=2p(z),\ p(\lambda z)=\lambda^{2}p(z),\label{eq:}\end{equation}
 \begin{equation}
p(z_{1}\pm z_{2})=p(z_{1})+p(z_{2})\pm z_{1}\cdot z_{2}\label{eq:}\end{equation}
\begin{equation}
p(z_{1}+z_{2})+p(z_{1}-z_{2})=2(p(z_{1})+p(z_{2})),\label{eq:}\end{equation}
for every $\lambda\in\mathbb{R}$, $z,z_{1},z_{2}\in Z$.

On $Z$ we fix a topology compatible with the natural duality such
as the strong$\times$weakly-star topology. In the sequel all topological
notions will be understood with respect to this fixed topology on
$Z$ if not otherwise specified.

For $\emptyset\neq A\subset Z$ let\begin{equation}
p_{A}=p+i_{A},\label{eq:}\end{equation}
where $i_{A}(z)=0$, if $z\in A$, $i_{A}(z)=\infty$, otherwise;
is the indicator of $A$.

The \emph{Fitzpatrick function} of $A$ is $h_{A}:Z\rightarrow\mathbb{R}\cup\{\infty\}$
defined by\[
h_{A}=p_{A}^{*}.\]

An alternative form for $h_{A}$ is given by\begin{equation}
h_{A}(z)=\sup_{a\in A}\{ z\cdot a-p(a)\}=p(z)-\inf_{a\in A}p(z-a),\ z\in Z.\label{ff}\end{equation}

The conjugate of $h_{A}$ \begin{equation}
\varphi_{A}=h_{A}^{*}={\rm cl}\ {\rm co}\ p_{A},\label{eq:}\end{equation}
is called the \emph{Penot function} of $A$ and it represents the
greatest proper convex lower semicontinuous function majorized by
$p_{A}$.

In the sequel, for $h:Z\rightarrow\mathbb{R}\cup\{\infty\}$, the
following set notation will be frequently used\[
\{ h=p\}=\{ z\in Z;\ h(z)=p(z)\},\]
\[
\{ h\ge p\}=\{ z\in Z;\ h(z)\ge p(z)\},\]
\[
\{ h\le p\}=\{ z\in Z;\ h(z)\le p(z)\}.\]

\strut

\noindent PROPOSITION 2.1. \emph{For every monotone $A\subset Z$
we have}

\emph{i) $\varphi_{A}(z)\ge p(z)$, for every $z\in Z$,}

\emph{ii) $A\subset\{\varphi_{A}=p\}$,}

\emph{iii) $A\subset\{ h_{A}=p\}$,}

\emph{iv) $D(A)\times X^{*}\subset\{ h_{A}\ge p\}$.}

\strut

\emph{Proof}. iv) If $z=(x,x^{*})\in D(A)\times X^{*}$ then there
exists $\alpha^{*}\in Ax$. Let $a=(x,\alpha^{*})\in A$. We have
$z\cdot a-p(a)=p(z)$. According to (\ref{ff}), this yields $h_{A}(z)\ge p(z)$,
i.e., $D(A)\times X^{*}\subset\{ h_{A}\ge p\}$.

iii) Let $z\in A$. Pick $a=z$ in (\ref{ff}) to find $h_{A}(z)\ge p(z)$.
Since $A$ is monotone we get $p(z-a)=p(a)-z\cdot a+p(z)\ge0$, for
every $a\in A$. This gives us $h_{A}(z)\le p(z)$. Therefore $A\subset\{ h_{A}=p\}$
and $h_{A}\le p_{A}$ in $Z$. Because $h_{A}$ is proper convex lower
semicontinuous this yields\begin{equation}
h_{A}\le\varphi_{A}\ {\rm in}\ Z,\label{eq:}\end{equation}
for every $A$ monotone.

Let ${\cal A}$ be a maximal monotone extension of $A$. By Theorem
3.3. below, we know that $h_{{\cal A}}\ge p$ in $Z$. We have $p_{A}\ge p_{{\cal A}}$
and\begin{equation}
\varphi_{A}\ge\varphi_{{\cal A}}\ge h_{{\cal A}}\ge p\ {\rm in}\ Z,\label{eq:}\end{equation}
i.e., i) holds.

Subpoint ii) is straight forward from i) and $\varphi_{A}\le p_{A}$
in $Z$. The proof is complete. \hfill $\square$

\strut

For other properties of $h_{A}$, $\varphi_{A}$ see {[}12, Proposition
2].

\strut

\strut

\noindent \textbf{3. Representability and maximality}

\strut

\noindent DEFINITION 3.1. A multi-valued operator $A$ is called \emph{representable}
in $Z=X\times X^{*}$ if there is a proper convex lower semicontinuous
$h:Z\rightarrow\mathbb{R}\cup\{\infty\}$ such that

i) $h(z)\ge p(z)$, for every $z\in Z$, i.e., $\{ h\ge p\}=Z$,

ii) $z\in A$ iff $h(z)=p(z)$, i.e., $A=\{ h=p\}$.

\smallskip

A function $h$ with properties i), ii) is called a \emph{representative}
of $A$. Notice that if $h$ is a representative of $A$ then from
$A\subset\{ h=p\}$ we get $h\le p_{A}$ in $Z$ followed by\begin{equation}
h\le\varphi_{A},\ h^{*}\ge h_{A}\ {\rm in}\ Z.\label{ra*>ha}\end{equation}

\strut

\noindent LEMMA 3.1. ({[}5, Proposition 4]) \textbf{}\emph{Every representable
operator $A$ is monotone.}

\strut

\emph{Proof}. Since $A$ is representable, there exists $h:Z\rightarrow\mathbb{R}\cup\{\infty\}$
such that $h\ge p$ in $Z$ and $z\in A$ iff $h(z)=p(z)$. Therefore,
from the convexity of $h$ we get that, for every $z_{1},z_{2}\in A=\{ h=p\}$\[
\frac{1}{2}p(z_{1})+\frac{1}{2}p(z_{2})-\frac{1}{4}p(z_{1}-z_{2})=p(\frac{1}{2}z_{1}+\frac{1}{2}z_{2})\]
\begin{equation}
\le h(\frac{1}{2}z_{1}+\frac{1}{2}z_{2})\le\frac{1}{2}h(z_{1})+\frac{1}{2}h(z_{2})=\frac{1}{2}p(z_{1})+\frac{1}{2}p(z_{2}),\label{eq:}\end{equation}
that is $p(z_{1}-z_{2})\ge0$, for every $z_{1},z_{2}\in A$.\hfill 
$\square$

\strut

We prove that an operator is representable iff its Penot function
becomes a representative.

\strut

\noindent THEOREM 3.2. $A$ \emph{is representable iff} $\varphi_{A}$
\emph{is a representative of} $A$.

\strut

\emph{Proof}. For the direct implication let $h$ be a representative
of $A$. From $A=\{ h=p\}$ we know that $h\le p_{A}$ in $Z$. Therefore
\begin{equation}
\varphi_{A}\ge h\ge p,\ {\rm in}\ Z.\label{eq:}\end{equation}
Combined with $A\subset\{\varphi_{A}=p\}$ the previous inequality
shows that $A=\{\varphi_{A}=p\}$, which implies that $\varphi_{A}$
is a representative of $A$. The converse implication is plain.\hfill 
$\square$

\strut

For different proofs of the previous result see e.g. {[}4,5].

\strut

The following characterization of maximality in terms of representability
is due to Fitzpatrick {[}2, Theorem 3.8]. For the sake of completeness
we provide the reader with a short proof.

\strut

\noindent THEOREM 3.3. \emph{A multi-valued operator $A$ is maximal
monotone iff $h_{A}$ is a representative of} $A$.

\strut

\emph{Proof}. If $A$ is maximal monotone then for every $z\not\in A$
there exists an $a\in A$ such that $p(z-a)<0$. Hence, from (\ref{ff})
we have $h_{A}(z)>p(z)$ for every $z\not\in A$. Since $h_{A}(z)=p(z)$,
for every $z\in A$ (see Proposition 2.1. iii)) this implies that
$h_{A}\ge p$ in $Z$ and $h_{A}(z)=p(z)$ iff $z\in A$, that is
$h_{A}$ is a representative of $A$.

Conversely, from Lemma 3.1. or from $h_{A}(z)=p(z)$ for every $z\in A$
and (\ref{ff}) we get $p(z-a)\ge0$ for every $z,a\in A$, i.e.,
$A$ is monotone.

Take $z_{0}\in Z$ such that $p(z_{0}-a)\ge0$ for every $a\in A$.
Again, from (\ref{ff}) we find $h_{A}(z_{0})=p(z_{0})$, that is
$z_{0}\in A$, since $h_{A}$ is a representative of $A$. We showed
that $A$ is maximal monotone. The proof is complete.\hfill  $\square$

\strut

Clearly, every maximal monotone operator is representable. The question
whether the converse holds appears naturally in this context. The
following characterization of maximality in terms of representability
appeared first in Voisei {[}11, Theorem 2.3]. For the sake of convenience
we provide the reader with a simpler proof.

\strut

\noindent THEOREM 3.4. $A$ \emph{is maximal monotone iff $A$ is
representable and} $h_{A}\ge p$ in $Z$.

\strut

\emph{Proof}. The direct implication is trivial since $h_{A}$ is
a representative of $A$.

Conversely, we know that $A$ is monotone since it is representable.
According to Proposition 2.1., we have $A\subset\{ h_{A}=p\}$. To
conclude that $h_{A}$ is a representative of $A$ and consequently
that $A$ is maximal monotone, it is enough to prove that $\{ h_{A}=p\}\subset A$.
Let $z\in\{ h_{A}=p\}$. Clearly, $z$ is a global minimum point of
$h_{A}-p$. Therefore\begin{equation}
0\in\partial(h_{A}-p)(z),\label{cr}\end{equation}
where {}``$\partial$'' denotes the Clarke-Rockafellar subdifferential.
Since $p$ is continuously G\^{a}teaux differentiable with $\partial(-p(z))=\{-z\}$
and $h_{A}$ is convex, relation (\ref{cr}) reduces to $z\in\partial h_{A}(z)$
which can be equivalently restated as\[
h_{A}(z)+\varphi_{A}(z)=2p(z).\]
This implies $z\in\{\varphi_{A}=p\}=A$ because $h_{A}(z)=p(z)$ and
$A$ is representable. The proof is complete.\hfill  $\square$

\strut

\emph{Remark 3.5.} Condition $h_{A}\ge p$ in $Z$ is sometimes referred
to as $A$ is of negative infimum type or NI in $X\times X^{*}$.
Hence the previous characterization theorem can be restated as

\begin{center}

\begin{tabular}{|c|}
\hline 
Maximal Monotone = Representable+NI\tabularnewline
\hline
\end{tabular}

\end{center}

This characterization of maximality is more versatile because most
of the times the representability of operators is easily checked.
Usually, the difficulty lies into proving that the operators are of
NI type.

\strut

\strut

\noindent \textbf{4. A simple proof of Heisler's result}

\strut

Previous to the papers {[}11,12] there are two note-worthy results
for the maximality of the sum in a non-reflexive Banach space setting;
the result of Heisler for full-space domain operators and the result
of Phelps \& Simons (see {[}6, Theorem 7.2]) for linear single-valued
operators. We provide a simpler proof of the Heisler result in order
to observe the usefulness of the Fitzpatrick function and mention
that in the linear multi-valued case the problem has been completely
solved (see {[}12,13] or Theorem 5.13. below).

Recall the Heisler result

\strut

\noindent THEOREM 4.1. ({[}10, Theorem 37.4]) \emph{Let $X$ be a
Banach space possibly non-reflexive. If $A,B$ are maximal monotone
in $X\times X^{*}$ with $D(A)=D(B)=X$ then $A+B$ is maximal monotone.}

\strut

The previous proof of this result relies on a topological characterization
of maximal monotone operators with full-space domain. Our argument
is based on the following two lemmas

\strut

\noindent LEMMA 4.2. \emph{Let $A$ be monotone with $D(A)=X$. Then
$A$ is of NI type in $Z=X\times X^{*}$, i.e.,\begin{equation}
h_{A}(z)\ge p(z),\ {\rm for\ every}\ z\in Z,\label{eq:}\end{equation}
and $\{ h_{A}=p\}$ is the only maximal monotone extension of $A$
in} $X\times X^{*}$.

\strut

\emph{Proof}. According to Proposition 2.1. iv), $\{ h_{A}=p\}=Z$,
that is, $A$ is (NI). For the second part notice that $\{ h_{A}=p\}$
is representable monotone and every maximal monotone extension ${\cal A}$
of $A$ satisfies\begin{equation}
{\cal A}\subset\{ h_{A}\le p\}=\{ h_{A}=p\}.\label{eq:}\end{equation}
Therefore ${\cal A}=\{ h_{A}=p\}$ and $\{ h_{A}=p\}$ is the unique
maximal monotone extension of $A$.\hfill  $\square$

\strut

\noindent LEMMA 4.3. \emph{Let $A$ be monotone with $D(A)=X$. Then
$A$ is maximal monotone iff $A$ has convex values and $A$ is closed
with respect to the strong$\times$weakly-star convergence of bounded
nets in $X\times X^{*}$ given by $(x_{\alpha},x_{\alpha}^{*})\twoheadrightarrow(x,x^{*})$
$\Leftrightarrow$ $x_{\alpha}\rightarrow x$, strongly in $X$, $x_{\alpha}^{*}\rightarrow x^{*}$,
weakly star in $X^{*}$, and $(x_{\alpha}^{*})_{\alpha}$ is bounded
in} $X^{*}$.

\strut

\emph{Proof}. The direct implication is clear since every maximal
monotone operator has convex values and is closed with respect to
{}``$\twoheadrightarrow$''.

For the converse it is enough to show that $\{ h_{A}=p\}\subset A$.

Since $A$ is closed with respect to {}``$\twoheadrightarrow$''
we prove first that $A$ has weakly-star closed values. Indeed, if
$(x_{\alpha}^{*})_{\alpha}\subset Ax$, $x\in X$, and $x_{\alpha}^{*}\rightarrow x^{*}$
weakly-star in $X^{*}$ then $(x_{\alpha}=x,x_{\alpha}^{*})\twoheadrightarrow(x,x^{*})$
because $Ax$ is bounded. Therefore, $x^{*}\in Ax$, i.e., $Ax$ is
weakly-star closed for every $x\in X$.

Let $z=(x_{0},x_{0}^{*})\in\{ h_{A}=p\}$, that is, for every $(a,a^{*})\in A$\begin{equation}
\langle x_{0}-a,x_{0}^{*}-a^{*}\rangle\ge0.\label{15}\end{equation}
Assume by contradiction that $x_{0}^{*}\not\in Ax_{0}$. By a separation
theorem we find $v_{0}\in X$, such that \begin{equation}
\langle v_{0},x_{0}^{*}\rangle>\sup_{x^{*}\in Ax_{0}}\langle v_{0},x^{*}\rangle.\label{st}\end{equation}
For $t>0$, denote by $a_{t}=x_{0}+tv_{0}$ and take $a_{t}^{*}\in Aa_{t}$
in (\ref{15}) to find\begin{equation}
\langle v_{0},x_{0}^{*}-a_{t}^{*}\rangle\le0.\label{h1}\end{equation}

Notice that for $t\downarrow0$, $a_{t}\rightarrow x_{0}$, strongly
in $X$. Because $A$ is locally bounded at $x_{0}$, $(a_{t}^{*})_{t}$
is bounded in $X^{*}$. Therefore, by the Alaoglu Theorem, at least
on a subnet, we have $(a_{t},a_{t}^{*})\twoheadrightarrow(x_{0},a_{0}^{*})\in A$.

Pass to limit in (\ref{h1}) with $t\downarrow0$ to get $\langle v_{0},a_{0}^{*}\rangle\ge\langle v_{0},x_{0}^{*}\rangle$
which contradicts (\ref{st}). We proved $z\in A$, that is $A=\{ h_{A}=p\}$.
The proof is complete. $\square$

\strut

\emph{Proof of Theorem 4.1.} It is straight forward to show that if
$A,B$ are closed with respect to {}``$\twoheadrightarrow$'' and
have convex values then $A+B$ is {}``$\twoheadrightarrow$'' closed
and has convex values, because $A,B$ are locally bounded. According
to Lemma 4.3., $A+B$ is maximal monotone.\hfill  \hfill $\square$

\strut

\strut

\noindent \textbf{5. Calculus rules for representable and maximal
monotone operators}

\strut

Our first concern in this section is the representability of $T:=L^{*}ML:X\rightrightarrows X^{*}$,
where $X,Y$ are Banach spaces, $L:X\rightarrow Y$ is linear bounded,
and $M:Y\rightrightarrows Y^{*}$ is representable. Let $r_{M}$ be
a representative of $M$.

Our choice for a representative of $T$ is $r:X\times X^{*}\rightarrow\overline{\mathbb{R}}$,\begin{equation}
r(x,x^{*})=\inf\{ r_{M}(Lx,y^{*});\ L^{*}y^{*}=x^{*}\},\ (x,x^{*})\in X\times X^{*}.\label{rcr}\end{equation}
Notice that $r\ge p$ in $Z$, $T\subset\{ r=p\}$, and $T=\{ r=p\}$
whenever the {}``inf'' in the definition of $r$ is attained for
all $(x,x^{*})\in D(r)$. Therefore, it is enough to study conditions
which assures that the {}``inf'' in (\ref{rcr}) becomes a {}``min''.

For a subset $S$ of a Banach space $X$ we denote by $^{i}S$ the
relative algebraic interior of $S$. We define $^{{\rm ic}}S={}^{i}S$
if the affine hull of $S$ is closed and $^{{\rm ic}}S=\emptyset$
otherwise. 

\strut

\noindent THEOREM 5.1. \emph{Let $X,Y$ be Banach spaces, $L:X\rightarrow Y$
be linear bounded, $M:Y\rightrightarrows Y^{*}$ be representable,
and $r_{M}$ be a representative of $M$. If}

\emph{\begin{equation}
0\in\,^{{\rm ic}}(R(L)-P_{Y}D(r_{M}^{*})),\label{eq:}\end{equation}
}

\noindent \emph{then $T:=L^{*}ML:X\rightrightarrows X^{*}$ is representable.
Here $P_{Y}:Y\times Y^{*}\rightarrow Y$, $P_{Y}(y,y^{*})=y$, $(y,y^{*})\in Y\times Y^{*}$
is the projection of $Y\times Y^{*}$ onto $Y$, $L^{*}:Y^{*}\rightarrow X^{*}$
denotes the adjoint of $L$, $R(L)$ is the range of $L$, and $r_{M}^{*}$
stands for the convex conjugate with respect to the natural duality.}

\strut

\emph{Proof}. Consider $\varphi:X\times X^{*}\rightarrow\overline{\mathbb{R}}$,\begin{equation}
\varphi(x,x^{*})=\inf\{ r_{M}^{*}(Lx,y^{*});\ L^{*}y^{*}=x^{*}\},\ (x,x^{*})\in X\times X^{*}.\label{eq:}\end{equation}
Since $\varphi(x,x^{*})=\inf\{ r_{M}^{*}(y,y^{*});\ (y,y^{*})\in C(x,x^{*})\}$,
where the process $C\subset X\times X^{*}\times Y\times Y^{*}$ is
defined by\begin{equation}
(x,x^{*},y,y^{*})\in C\ {\rm iff}\ y=Lx,\ x^{*}=L^{*}y^{*},\label{defC}\end{equation}
$R(C)=R(L)\times Y^{*}$, and condition $0\in\,^{{\rm ic}}(R(L)-P_{Y}D(r_{M}^{*}))$
is equivalent to $0\in\,^{{\rm ic}}(R(C)-D(r_{M}^{*}))$, according
to {[}15, T 2.8.6], we may apply the chain rule to get\begin{equation}
\varphi_{Z^{*}}^{*}(x^{*},x^{**})=\min\{ r_{M}^{**}(y^{*},y^{**});\ (x^{*},x^{**})\in C^{*}(y^{*},y^{**})\},\ (x^{*},x^{**})\in X^{*}\times X^{**}.\label{me0}\end{equation}
Here $\varphi_{Z^{*}}^{*}$ denotes the convex conjugate of $\varphi$
in $Z^{*}=X^{*}\times X^{**}$ and is weakly-star lower semicontinuous
in $Z^{*}$ which makes $\varphi^{*}=\varphi_{Z^{*}}^{*}/_{Z}$ weakly$\times$weakly-star
lower semi-continuous in $Z$.

The adjoint of $C$ is given by \begin{equation}
(y^{*},y^{**},x^{*},x^{**})\in C^{*}\ {\rm iff}\ y^{**}=L^{**}x^{**},\ x^{*}=L^{*}y^{*}.\label{defC*}\end{equation}
By the bi-conjugate formula for $x^{**}=x\in X$ relation (\ref{me0})
becomes\begin{equation}
r(x,x^{*})=\varphi^{*}(x,x^{*})=\min\{ r_{M}(Lx,y^{*});\ L^{*}y^{*}=x^{*}\},\ (x,x^{*})\in X\times X^{*}.\label{me}\end{equation}
Relation (\ref{me}) shows that $r$ is a representative of $T$ and
consequently $T$ is representable.\hfill  $\square$

\strut

\noindent PROPOSITION 5.2. \emph{Let $X$ be a Banach space, $L:X\rightarrow X\times X$,
$Lx=(x,x)$, $x\in X$, and $U,V$ be convex subsets of $X$. Then}\begin{equation}
0\in\,^{i}(R(L)-U\times V)\Leftrightarrow0\in\,^{i}(U-V).\label{28}\end{equation}
 \begin{equation}
0\in\,^{{\rm ic}}(R(L)-U\times V)\Leftrightarrow0\in\,^{{\rm ic}}(U-V).\label{echiv ic}\end{equation}

\strut

\emph{Proof}. Consider the {}``difference function'' ${\cal D}:X\times X\rightarrow X$,
\[
{\cal D}(x_{1},x_{2})=x_{2}-x_{1},\ x_{1},x_{2}\in X.\]
Notice that $\mbox{Ker}{\cal D}=R(L)$, and let {}``aff'' denote
the affine hull of a subset in $X$ or $X\times X$. We have\begin{equation}
F:=\mbox{aff}(R(L)-U\times V)=R(L)-\mbox{aff}U\times\mbox{aff}V,\label{eq:}\end{equation}
\begin{equation}
{\cal D}(F)=\mbox{aff}(U-V).\label{eq:}\end{equation}
According to {[}15, Corollary 1.3.15] ${\cal D}(F)$ is closed iff
$F=F+\mbox{Ker}{\cal D}$ is closed. Therefore, $\mbox{aff}(R(L)-U\times V)$
is closed iff $\mbox{aff}(U-V)$ is closed.

Since $U,V$ are convex, condition $0\in\,^{i}(R(L)-U\times V)$ is
equivalent to $\bigcup\limits _{n\ge1}n(R(L)-U\times V)$ is a linear
subspace and $0\in\,^{i}(U-V)$ is equivalent to $\bigcup\limits _{n\ge1}n(U-V)$
is a linear subspace (see {[}15, (1.1)]).

But\begin{equation}
(x_{1},x_{2})\in R(L)-U\times V\ {\rm iff}\ {\cal D}(x_{1},x_{2})=x_{2}-x_{1}\in U-V,\label{eq:}\end{equation}
which shows that $\bigcup\limits _{n\ge1}n(R(L)-U\times V)$ is a
linear subspace iff $\bigcup\limits _{n\ge1}n(U-V)$ is a linear subspace.
Hence $0\in\,^{i}(R(L)-U\times V)$ iff $0\in\,^{i}(U-V)$. The proof
is complete.\hfill  $\square$

\strut

For a generalization of (\ref{28}) see {[}16, Proposition 2.1].

\eject

\noindent THEOREM 5.3. \emph{Let $X$ be a Banach space and $A,B:X\rightrightarrows X^{*}$
be representable with}

\emph{\begin{equation}
0\in\,^{{\rm ic}}(P_{X}D(r_{A}^{*})-P_{X}D(r_{B}^{*})),\label{me1}\end{equation}
where $r_{A},r_{B}$ are representatives of $A,B$ and $P_{X}$ stands
for the projection of $X\times X^{*}$ onto $X$. Then $A+B$ is representable.}

\strut

\emph{Proof}. First argument. We apply Theorem 5.1. for $X$, $Y=X\times X$,
$Lx=(x,x)$, $x\in X$, and $M(x_{1},x_{2})=Ax_{1}\times Bx_{2}$,
$(x_{1},x_{2})\in D(M)=D(A)\times D(B)$ for which $L^{*}ML=A+B$,\begin{equation}
r_{M}(x_{1},x_{2},x_{1}^{*},x_{2}^{*})=r_{A}(x_{1},x_{1}^{*})+r_{B}(x_{2},x_{2}^{*}),\ x_{1},x_{2}\in X,\ x_{1}^{*},x_{2}^{*}\in X^{*},\label{eq:}\end{equation}
is a representative of $M$ with\begin{equation}
r_{M}^{*}(x_{1},x_{2},x_{1}^{*},x_{2}^{*})=r_{A}^{*}(x_{1},x_{1}^{*})+r_{B}^{*}(x_{2},x_{2}^{*}),\ x_{1},x_{2}\in X,\ x_{1}^{*},x_{2}^{*}\in X^{*},\label{eq:}\end{equation}
$P_{Y}D(r_{M}^{*})=P_{X}D(r_{A}^{*})\times P_{X}D(r_{B}^{*})$, and
according to Proposition 5.2., condition $0\in\,^{{\rm ic}}(P_{X}D(r_{A}^{*})-P_{X}D(r_{B}^{*}))$
is equivalent to $0\in\,^{{\rm ic}}(R(L)-P_{Y}D(r_{M}^{*}))$.

Second argument. Let $(x_{0},x_{0}^{*})\in\{\varphi_{A+B}=p\}$. Since
$\varphi_{A+B}=h_{A+B}^{*}$ we have that for every $(u,u^{*})\in Z$\begin{equation}
h_{A+B}(u,u^{*})-\langle u,x_{0}^{*}\rangle-\langle x_{0},u^{*}\rangle+\langle x_{0},x_{0}^{*}\rangle\ge0.\label{l1}\end{equation}
Let ${\cal X}=X\times X^{*}\times X^{*}$, ${\cal Y}=X$ and consider
the function $\Phi:{\cal X}\times{\cal Y}\rightarrow\mathbb{R}\cup\{\infty\}$
given by\[
\Phi(x,x^{*},z^{*};y)=r_{A}^{*}(x+y,x^{*})+r_{B}^{*}(x,z^{*})-\langle x,x_{0}^{*}\rangle-\langle x_{0},x^{*}+z^{*}\rangle+\langle x_{0},x_{0}^{*}\rangle,\]
$x,y\in X$, $x^{*},z^{*}\in X^{*}$.

Notice that, since $r_{A}^{*}\ge h_{A}$, $r_{B}^{*}\ge h_{B}$ (see
(\ref{ra*>ha})), we get

\[
\Phi(x,x^{*},z^{*};0)\ge h_{A}(x,x^{*})+h_{B}(x,z^{*})-\langle x,x_{0}^{*}\rangle-\langle x_{0},x^{*}+z^{*}\rangle+\langle x_{0},x_{0}^{*}\rangle\]
\[
\ge h_{A+B}(x,x^{*}+z^{*})-\langle x,x_{0}^{*}\rangle-\langle x_{0},x^{*}+z^{*}\rangle+\langle x_{0},x_{0}^{*}\rangle\ge0,\]
for every $(x,x^{*},z^{*})\in{\cal X}$, i.e., $\inf\limits _{\chi\in{\cal X}}\Phi(\chi,0)\ge0$.

If $P_{{\cal Y}}(\chi,y)=y$, $(\chi,y)\in{\cal X}\times{\cal Y}$,
is the projection of ${\cal X}\times{\cal Y}$ onto ${\cal Y}$ then
$P_{{\cal Y}}D(\Phi)=P_{X}D(r_{A}^{*})-P_{X}D(r_{B}^{*})$ and condition
(\ref{me1}) spells $0\in\,^{{\rm ic}}(P_{{\cal Y}}D(\Phi))$. 

This allows us to apply the fundamental duality formula (see e.g.
{[}15, Theorem 2.7.1 (vii)]) to get\begin{equation}
\inf\limits _{\chi\in{\cal X}}\Phi(\chi,0)=\max_{y^{*}\in Y^{*}=X^{*}}(-\Phi^{*}(0,y^{*}))\ge0.\label{eq:}\end{equation}
Therefore, there exists $y^{*}\in X^{*}$ such that\begin{equation}
\Phi^{*}(0,y^{*})=\sup\{\langle y,y^{*}\rangle-\Phi(x,x^{*},z^{*};y);\ x,y\in X,\ x^{*},z^{*}\in X^{*}\}\le0,\label{eq:}\end{equation}
that is \begin{equation}
r_{A}^{*}(x+y,x^{*})+r_{B}^{*}(x,z^{*})-\langle x,x_{0}^{*}\rangle-\langle x_{0},x^{*}+z^{*}\rangle+\langle x_{0},x_{0}^{*}\rangle-\langle y,y^{*}\rangle\ge0,\label{eq:}\end{equation}
for every $x,y\in X$, $x^{*},z^{*}\in X^{*}$.

Using the substitution $x+y=z$, we find\begin{equation}
r_{A}^{*}(z,x^{*})+r_{B}^{*}(x,z^{*})-\langle x,x_{0}^{*}\rangle-\langle x_{0},x^{*}+z^{*}\rangle+\langle x_{0},x_{0}^{*}\rangle-\langle z-x,y^{*}\rangle\ge0,\label{eq:}\end{equation}
for every $x,y\in X$, $x^{*},z^{*}\in X^{*}$, or\[
[\langle z,y^{*}\rangle+\langle x_{0},x^{*}\rangle-r_{A}^{*}(z,x^{*})]+[\langle x,x_{0}^{*}-y^{*}\rangle+\langle x_{0},z^{*}\rangle-r_{B}^{*}(x,z^{*})]\le\langle x_{0},x_{0}^{*}\rangle,\]
for every $x,y\in X$, $x^{*},z^{*}\in X^{*}$, that is\begin{equation}
r_{A}(x_{0},y^{*})+r_{B}(x_{0},x_{0}^{*}-y^{*})\le\langle x_{0},x_{0}^{*}\rangle.\label{imp}\end{equation}
Because $r_{A},r_{B}$ are representatives of $A,B$ relation (\ref{imp})
is equivalent to\begin{equation}
r_{A}(x_{0},y^{*})=\langle x_{0},y^{*}\rangle,\ r_{B}(x_{0},x_{0}^{*}-y^{*})=\langle x_{0},x_{0}^{*}-y^{*}\rangle,\label{eq:}\end{equation}
that is, $(x_{0},y^{*})\in A$, $(x_{0},x_{0}^{*}-y^{*})\in B$. Hence
$(x_{0},x_{0}^{*})\in A+B.$ We proved that $\{\varphi_{A+B}=p\}\subset A+B$
and this is enough in order to conclude that $A+B$ is representable.\hfill 
$\square$

\strut

\emph{Remark 5.4.} The typical example of a representative of $A$
is provided by the Penot function $\varphi_{A}$. Therefore, in a
particular case, Theorems 5.1, 5.3 can be restated as 

\strut

\noindent COROLLARY 5.5. \emph{Let $X,Y$ be two Banach spaces, $L:X\rightarrow Y$
be linear bounded, and $M:Y\rightrightarrows Y^{*}$ be representable.
If\begin{equation}
0\in\,^{{\rm ic}}(R(L)-P_{Y}D(h_{M}))\label{c1}\end{equation}
then $T:=L^{*}ML:X\rightrightarrows X^{*}$ is representable.}

\strut

\noindent COROLLARY 5.6. \emph{Let $X$ be a Banach space and $A,B:X\rightrightarrows X^{*}$
be representable with\[
0\in\,^{{\rm ic}}(P_{X}D(h_{A})-P_{X}D(h_{B})).\]
Then $A+B$ is representable.}

\strut

\noindent PROPOSITION 5.7. \emph{If $M:Y\rightrightarrows Y^{*}$
is monotone in the Banach space $Y$, $D(M)$ is closed convex, and\begin{equation}
M=M+N_{D(M)}\label{e1}\end{equation}
then $M$ is of NI type and $D(M)=P_{Y}D(h_{M})$. Here $N_{D(M)}$
stands for the convex normal cone to $D(M)$.}

\strut

\emph{Proof}. For every $y\in D(M)$ there is $y^{*}\in Y^{*}$ such
that $(y,y^{*})\in M\subset\{ h_{M}=p\}\subset D(h_{M})$, that is,
$D(M)\subset P_{Y}D(h_{M})$.

Conversely, let $y\in P_{Y}D(h_{M})$, that is, $h_{M}(y,y^{*})<\infty$,
for some $y^{*}\in Y^{*}$. Hence, for every $(m,m^{*})\in M$, we
have\begin{equation}
\langle y-m,m^{*}\rangle+\langle m,y^{*}\rangle\le C<\infty.\label{e2}\end{equation}
From (\ref{e1}), (\ref{e2}) and because $N_{D(M)}(y)$ is a cone
for every $y\in D(M)$, we get\begin{equation}
t\langle y-m,n^{*}\rangle+\langle y-m,m^{*}\rangle+\langle m,y^{*}\rangle\le C<\infty,\label{e3}\end{equation}
for every $t>0$, $m\in D(M)$, $m^{*}\in Mm$, $n^{*}\in N_{D(M)}(m)$.

From (\ref{e3}) it yields that $\langle y-m,n^{*}\rangle\le0$, for
every $m\in D(M)$, $n^{*}\in N_{D(M)}(m)$, i.e., $(y,0)$ is monotonically
related to the graph of the maximal monotone operator $N_{D(M)}$.
Therefore, $(y,0)\in N_{D(M)}$, that is, $y\in D(M)$. We proved
$P_{Y}D(h_{M})=D(M)$, i.e., $D(h_{M})\subset D(M)\times X^{*}$.
According to Proposition 2.1. iv) this implies $D(h_{M})\subset\{ h_{M}\ge p\}$,
that is $M$ is NI.\hfill  $\square$

\strut

\emph{Remark 5.8.} Condition (\ref{e1}) is satisfied whenever $M$
is maximal monotone, since $N_{D(M)}$ is monotone, $0\in N_{D(M)}(y)$
for every $y\in D(M)$, and $M\subset M+N_{D(M)}$. Therefore, every
maximal monotone $M$ with $D(M)$ closed convex has $P_{Y}D(h_{M})=D(M)$.

\strut

\noindent THEOREM 5.9. \emph{Let $X,Y$ be Banach spaces, $L:X\rightarrow Y$
be linear bounded, $M:Y\rightrightarrows Y^{*}$ be maximal monotone,
and $T:=L^{*}ML:X\rightrightarrows X^{*}$.}

\emph{($\alpha$) If $D(T)$ is closed convex and \begin{equation}
N_{D(T)}=L^{*}N_{D(M)}L,\label{e4}\end{equation}
then $T$ is of NI type.}

\emph{($\beta$) If $D(M)$ is closed convex and\begin{equation}
0\in\,^{{\rm ic}}(R(L)-D(M)),\label{e7}\end{equation}
}

\noindent \emph{then $T$ is maximal monotone.}

\emph{($\gamma$) If $D(T)$ is closed and $R(L)\cap{\rm int}D(M)\neq\emptyset$
then $T$ is maximal monotone.}

\strut

\emph{Proof}. ($\alpha$) Since $M$ is maximal monotone we know that
$M=M+N_{D(M)}$. We find $M(Lx)=M(Lx)+N_{D(M)}(Lx)$ and\[
Tx=L^{*}M(Lx)=L^{*}M(Lx)+L^{*}N_{D(M)}(Lx)=Tx+N_{D(T)}(x),\]
for every $x\in D(T)=L^{-1}(D(M))$, that is, $T=T+N_{D(T)}$ and,
according to Proposition 5.7., $T$ is of NI type.

($\beta$) Since $D(M)$ is closed convex, $D(T)=L^{-1}(D(M))$ is
closed convex, $D(M)=P_{Y}D(h_{M})$, and (\ref{e7}) becomes (\ref{c1}),
and so, by Corollary 5.5., $T$ is representable. Also, \begin{equation}
i_{D(T)}(x)=\inf\{ i_{D(M)}(y);\ y=Lx\},\ x\in D(T).\label{t86}\end{equation}
Taking into account (\ref{e7}), we may apply the chain rule {[}15,
Theorem 2.8.6 (v)] to get that\begin{equation}
N_{D(T)}=L^{*}N_{D(M)}L.\label{eq:}\end{equation}
According to ($\alpha$), $T$ is of NI type. Hence $T$ is maximal
monotone.

($\gamma$) Because $M$ is maximal monotone with $\mbox{int}D(M)\neq\emptyset$,
$\mbox{int}D(M)$, $\overline{D(M)}$ are convex, $\mbox{int}D(M)=\mbox{int}\overline{D(M)}$,
$\overline{D(M)}=\overline{\mbox{int}D(M)}$ (see e.g {[}10, Theorem
18.4]), $\mbox{int}D(M)=\mbox{int}P_{Y}D(h_{M})$ (see e.g. {[}9,
Theorem 2.2.]), $R(L)-P_{Y}D(h_{M})$ contains $0$ in its interior,
and (\ref{c1}) follows making $T$ representable.

We prove that \begin{equation}
D(T)=L^{-1}(\overline{D(M)}).\label{e10}\end{equation}
The direct inclusion is plain since $L$ is continuous and $D(T)$
is closed.

Conversely, let $x_{0}\in L^{-1}(\overline{D(M)})$, that is, $Lx_{0}\in\overline{D(M)}$.
Without loss of generality we may assume that $0\in\mbox{int}D(M)$
and $0\in M0$. Then $\lambda Lx_{0}\in D(M)$, for every $0\le\lambda<1$
(see e.g. {[}15, Theorem 1.1.2]), i.e., $\lambda x_{0}\in D(T)$,
for $0\le\lambda<1$. Letting $\lambda\uparrow1$, we find $x_{0}\in\overline{D(T)}=D(T)$.

Relation (\ref{e10}) shows that $D(T)$ is closed convex.

Again, from the chain rule {[}15, Theorem 2.8.6 (iii)] applied for

\[
i_{D(T)}(x)=i_{L^{-1}(\overline{D(M)})}(x)=\mbox{inf}\{ i_{\overline{D(M)}}(y);\ y=Lx\},\ x\in X.\]
we get (\ref{e4}), that is $T$ is NI and this is sufficient in order
to conclude.\hfill  $\square$ 

\strut

\noindent THEOREM 5.10. \emph{Let $A,B$ be maximal monotone operators
in the Banach space $X$.}

\emph{($\alpha$) If $D(A)\cap D(B)$ is closed convex and \begin{equation}
N_{D(A)\cap D(B)}=N_{D(A)}+N_{D(B)},\label{e5}\end{equation}
then $A+B$ is of NI type.}

\emph{($\beta$) If $D(A)$, $D(B)$ are closed convex and\begin{equation}
0\in\,^{{\rm ic}}(D(A)-D(B)),\label{eq:}\end{equation}
then $A+B$ is maximal monotone,}

\emph{($\gamma$) If $D(A)\cap D(B)$ is closed, $\overline{D(A)}$
is convex, and $D(A)\cap{\rm int}D(B)\neq\emptyset$ then $\overline{D(A)}\cap\overline{D(B)}=D(A)\cap D(B)$
and $A+B$ is maximal monotone.}

\emph{($\delta$) If $D(A)$ is closed convex and $D(A)\subset D(B)$
then $A+B$ is of NI type.}

\emph{($\epsilon$) If $D(A)$ is closed convex, $D(A)\subset D(B)$,
and $0\in\,^{{\rm ic}}(D(A)-P_{X}D(h_{B}))$ then $A+B$ is maximal
monotone.}

\strut

\emph{Proof}. Sub-points ($\alpha$), ($\beta$) are direct consequences
of Theorem 5.9. ($\alpha$), ($\beta$) applied for $Y=X\times X$,
$Lx=(x,x)$, $x\in X$, $L^{*}:Y^{*}=X^{*}\times X^{*}\rightarrow X^{*}$,
$L^{*}(x^{*},y^{*})=x^{*}+y^{*}$, $x^{*},y^{*}\in X^{*}$, $M(x_{1},x_{2})=Ax_{1}\times Bx_{2}$,
$(x_{1},x_{2})\in D(M)=D(A)\times D(B)$, for which $L^{*}ML=A+B$.
More precisely, subpoint ($\beta$) follows from Theorem 5.9. ($\beta$)
since $0\in\,^{{\rm ic}}(D(A)-D(B))$ iff $0\in\,^{{\rm ic}}(R(L)-D(M))$
(see (\ref{echiv ic})). For an alternative proof of ($\beta$) see
{[}12, Theorem 2].

($\gamma$) Without loss of generality assume that $0\in D(A)\cap{\rm int}D(B)$.
If $x\in\overline{D(A)}\cap\overline{D(B)}$, then, for every $0\le\lambda<1$,
$\lambda x\in\overline{D(A)}\cap{\rm int}D(B)\subset\overline{D(A)\cap D(B)}$.
Let $\lambda\uparrow1$ to find $x\in\overline{D(A)\cap D(B)}=D(A)\cap D(B)$,
that is \emph{$\overline{D(A)}\cap\overline{D(B)}=D(A)\cap D(B)$}
and consequently $D(A)\cap D(B)$ is convex\emph{.}

Therefore, for every $x\in D(A)\cap D(B)$\[
N_{D(A)\cap D(B)}(x)=N_{\overline{D(A)}\cap\overline{D(B)}}(x)\]
\[
=N_{\overline{D(A)}}(x)+N_{\overline{D(B)}}(x)=N_{D(A)}(x)+N_{D(B)}(x),\]
 i.e., (\ref{e5}) holds. The NI type follows from ($\alpha$) while
the representability is a consequence of \emph{$D(A)\cap{\rm int}D(B)\neq\emptyset$}
and Corollary 5.6.

($\delta$) Clearly, $D(A)\cap D(B)=D(A)$ is closed convex and since
$N_{D(A)}$ is maximal monotone we get \[
N_{D(A)\cap D(B)}=N_{D(A)}=N_{D(A)}+N_{\overline{co}D(B)}=N_{D(A)}+N_{D(B)},\]
i.e., according to ($\alpha$), $A+B$ is NI. Here {}``$\overline{co}$''
stands for the closed convex hull\emph{.}

($\epsilon$) Condition \emph{$0\in\,^{{\rm ic}}(D(A)-P_{X}D(h_{B}))$}
implies the representability of $A+B$. From ($\delta$) we know that
$A+B$ is NI, therefore $A+B$ is maximal monotone.\hfill  $\square$

\strut

\emph{Remark 5.11.} A recent results of Groh {[}3, Theorem 1.6] is
a particular case of our subpoint ($\epsilon$), for $A$ being a
subdifferential and $B$ having a non-empty interior.

\strut

The following result of Bauschke presents a different perspective
on the subject.

\strut

\noindent THEOREM 5.12. ({[}10, Theorem 39.1]) \emph{Let $A$ be maximal
monotone in the Banach space $X$ and $B:X\rightarrow X^{*}$ be linear
with $\langle Bx,x\rangle=0$, for every $x\in X$. Then $A+B$ is
maximal monotone.}

\strut

\emph{Proof}. It is easily checked that for every $(x,x^{*})\in X\times X^{*}$\begin{equation}
h_{A+B}(x,x^{*})=h_{A}(x,x^{*}+B^{*}x)=h_{A}(x,x^{*}-Bx),\label{eq:}\end{equation}
where $B^{*}=-B$ stands for the adjoint of $B$. This equality suffices
in order to conclude that $h_{A+B}$ is a representative of $A+B$
and $A+B$ is maximal monotone.\hfill  $\square$

\strut

Notice that under the assumptions of Bauschke's result we have\begin{equation}
h_{A+B}(x,x^{*})=\mbox{inf}\{ h_{A}(x,y^{*})+h_{B}(x,z^{*});\ y^{*}+z^{*}=x^{*}\}=(h_{A}\square_{2}h_{B})(x,x^{*}),\label{e14}\end{equation}
for every $(x,x^{*})\in X\times X^{*}$, since $h_{B}(x,z^{*})=0$,
iff $z^{*}=-B^{*}x$, $h_{B}(x,z^{*})=+\infty$, otherwise; where
{}``$\square_{2}$'' denotes the infimal convolution with respect
to the second variable.

It is worth noticing that equality (\ref{e14}) assures that $A+B$
is NI and that $A+B$ is maximal monotone whenever the infimal convolution
in (\ref{e14}) is exact. Unfortunately, in general (\ref{e14}) does
not hold even under the assumptions $D(A)=D(B)=X$ and $X$ is a Hilbert
space (see e.g. {[}7, Example 1]). Other cases in which an equality
of type (\ref{e14}) holds are given in the following theorem.

\strut

\noindent THEOREM 5.13. \emph{Let $X,Y$ be Banach spaces.}

\emph{($\alpha$) If $L:X\rightarrow Y$ is linear bounded, and $M:Y\rightrightarrows Y^{*}$
is maximal monotone with $\mbox{Graph}(M)$ convex in $X\times X^{*}$
and\begin{equation}
0\in\,^{{\rm ic}}(R(L)-D(M)),\label{e15}\end{equation}
}

\noindent \emph{then $T:=L^{*}ML:X\rightrightarrows X^{*}$ is maximal
monotone.}

\emph{($\beta$) If $A,B$ are maximal monotone operators in $X$
with $\mbox{Graph}(A)$, $\mbox{Graph}(B)$ convex and \begin{equation}
0\in\,^{{\rm ic}}(D(A)-D(B)),\label{e16}\end{equation}
then $A+B$ is maximal monotone.}

\strut

\emph{Proof}. ($\alpha$) We have\[
p_{T}(x,x^{*})=\min\{ p_{M}(Lx,y^{*});\ L^{*}y^{*}=x^{*}\}\qquad\qquad\qquad\]
\begin{equation}
\qquad=\min\{ p_{M}(y,y^{*});\ (y,y^{*})\in C(x,x^{*})\},\ (x,x^{*})\in X\times X^{*},\label{eq:}\end{equation}
where $C\subset X\times X^{*}\times Y\times Y^{*}$ is defined in
(\ref{defC}) with adjoint $C^{*}$ given by (\ref{defC*}).

Notice that $R(C)=R(L)\times Y^{*}$, $D(p_{M})=M$, $R(C)-D(p_{M})=(R(L)-D(M))\times Y^{*}$
and condition $0\in\,^{{\rm ic}}(R(C)-D(p_{M}))$ is equivalent to
(\ref{e15}). Moreover, $\mbox{Graph}M$ strongly closed and convex
in $X\times X^{*}$ makes $p_{M}$ proper convex strongly lower semicontinuous
in $X\times X^{*}$.

We apply the chain rule {[}15, Theorem 2.8.6 (v)] to get\begin{equation}
p_{T}^{*}(x^{*},x^{**})=\mbox{min}\{ p_{M}^{*}(y^{*},y^{**});\ (x^{*},x^{*})\in C^{*}(y^{*},y^{**})\},\label{e19}\end{equation}
$(x,x^{*})\in X\times X^{*}$. For $x^{**}=x\in X$ we find\begin{equation}
h_{T}(x,x^{*})=\mbox{min}\{ h_{M}(Lx,y^{*});\ L^{*}y^{*}=x^{*}\},\label{e20}\end{equation}
which implies that $h_{T}$ is a representative of $T$, i.e., $T$
is maximal monotone.

($\beta$) Again, take $Y=X\times X$, $Lx=(x,x)$, $x\in X$, $L^{*}:Y^{*}=X^{*}\times X^{*}\rightarrow X^{*}$,
$L^{*}(x^{*},y^{*})=x^{*}+y^{*}$, $x^{*},y^{*}\in X^{*}$, and $M(x_{1},x_{2})=Ax_{1}\times Bx_{2}$,
$(x_{1},x_{2})\in D(M)=D(A)\times D(B)$ or $\mbox{Graph}M=\mbox{Graph}A\times\mbox{Graph}B$.

Then $L^{*}ML=A+B$ is maximal monotone by the conclusion of ($\alpha$),
taking into consideration that, in this case, (\ref{e16}) is equivalent
to (\ref{e15}).\hfill  $\square$

\strut

\noindent COROLLARY 5.14. \emph{Let $X,Y$ be Banach spaces.}

\emph{($\alpha$) If $L:X\rightarrow Y$ is linear bounded, and $M:Y\rightrightarrows Y^{*}$
is linear maximal monotone with $R(L)-D(M)$ closed in $Y$, then
$T:=L^{*}ML:X\rightrightarrows X^{*}$ is maximal monotone.}

\emph{($\beta$) If $A,B$ are linear maximal monotone with $D(A)-D(B)$
closed in $X$ then $A+B$ is maximal monotone.}

\strut

\emph{Proof}. Condition (\ref{e15}) is equivalent to $R(L)-D(M)$
closed in $Y$, since $R(L)-D(M)$ is a subspace. Similarly, (\ref{e16})
becomes $D(A)-D(B)$ is closed in $X$. For a different proof of ($\beta$)
see {[}13].\hfill  $\square$

\strut

It is worth mentioning that in the linear case the qualification constraints
contained in ($\alpha$), ($\beta$) cannot be further relaxed (see
e.g. {[}10] for a counter-example).

\strut

\strut

\strut

\noindent \textbf{Acknowledgments}

\strut

\noindent The author would like to thank Dr. C. Z\u{a}linescu for
several interesting comments made on a preliminary version of this
paper.

\strut

\strut

\noindent \textbf{References}

\begin{enumerate}
\item Borwein, J.: Maximal monotonicity via convex analysis, \emph{Journal
of Convex Analysis} \textbf{13/14} (2006).
\item Fitzpatrick, S.: Representing monotone operators by convex functions,
Workshop/Miniconference on Functional Analysis and Optimization (Canberra,
1988), in: \emph{Proc. Centre Math. Anal. Austral. Nat. Univ.} \textbf{20},
Austral. Nat. Univ., Canberra, 1988, pp. 59--65. 
\item Groh K.: On monotone operators and forms, \emph{Journal of Convex
Analysis} \textbf{12}(2) (2005), 417--429.
\item Martinez-Legaz, J.E. and Svaiter, B.F.: Monotone operators representable
by l.s.c. convex functions, \emph{Set-Valued Analysis} \textbf{13}
(2005), 21--46.
\item Penot, J.-P,: The relevance of convex analysis for the study of monotonicity,
(English. English summary) \emph{Nonlinear Anal.} \textbf{58}(7-8)
(2004), 855--871.
\item Phelps, R.R. and Simons S.: Unbounded linear monotone operators on
nonreflexive Banach spaces (English. English summary), \emph{Journal
of Convex Analysis} \textbf{5}(2) (1998), 303--328.
\item Penot, J.-P and Z\u alinescu, C.: Some problems about the representation
of monotone operators by convex functions, \emph{Anziam J.} \textbf{47}(1)
(2005), 1--20.
\item Rockafellar, R.T.: On the maximality of sums of nonlinear monotone
operators, \emph{Trans. Amer. Math. Soc.} \textbf{149} (1970), 75--88.
\item Simons, S.: Dualized and scaled Fitzpatrick functions, \emph{Proc.
A.M.S.} \textbf{134}(10) (2006), 2983--2987.
\item Simons, S.: \emph{Minimax and monotonicity}, Lecture Notes in Mathematics,
1693. Springer-Verlag, Berlin, 1998.
\item Voisei, M.D.: A maximality theorem for the sum of maximal monotone
operators in non-reflexive Banach Spaces, \emph{Math. Sci. Res. J.},
\textbf{10}(2) (2006), 36--41.
\item Voisei, M.D.: Calculus rules for maximal monotone operators in general
Banach spaces, preprint 2006.
\item Voisei, M.D.: The sum theorem for linear maximal monotone operators,
\emph{Math. Sci. Res. J.} \textbf{10}(4) (2006), 83--85.
\item Z\u alinescu C.: A new proof of the maximal monotonicity of the sum
using the Fitzpatrick function, Variational analysis and applications,
in: Nonconvex Optim. Appl., \textbf{79}, Springer, New York, 2005,
pp. 1159--1172.
\item Z\u alinescu C.: \emph{Convex analysis in general vector spaces},
World Scientific Publishing Co., Inc., River Edge, NJ, 2002.
\item Z\u alinescu C.: On convex sets in general position, \emph{Linear
Algebra Appl.} \textbf{64} (1985), 191--198.
\end{enumerate}

\end{document}